\documentclass[10pt]{article}
\usepackage{cite,amsfonts,a4wide}
\usepackage{epsfig}
\usepackage[linguistics]{forest}
\usepackage{amsmath}
\usepackage{geometry}
\input amssym.def
\input amssym.tex

\newcommand {\N} {\mbox{$N$}}

\newcommand {\be}{\begin{equation}}
\newcommand {\e}{\end{equation}}
\newcommand {\bes}{\begin{displaymath}}
\newcommand {\es}{\end{displaymath}}
\newcommand {\bea}{\begin{align}}
\newcommand {\ea}{\end{align}}
\newcommand {\nit}{\noindent}

\newcommand {\bit}{\bibitem}
\newcommand {\kpr}{k^\prime}

\newcommand {\fract}[2]{\mbox{${\displaystyle{\frac{#1}{#2}}}$}}

\begin{document}


\begin{center}
{\bf \Large Pi Formulas: some smooth stones \\
on the beach of rough numbers\\
\vspace{0.5cm}
A.J. Macfarlane$\;$}\footnote{{\it e-mail}: 
a.j.macfarlane@damtp.cam.ac.uk}\\
\vskip .5cm
\begin{sl}
Centre for Mathematical Sciences, D.A.M.T.P.\\
Wilberforce Road, Cambridge CB3 0WA, UK
\vskip .5cm
\end{sl}
\end{center}

\begin{abstract}

This article is about Pi Formulas, infinite series of fractions which sum to
simple multiples of Pi. Each such one can be associated with  a unique set
$S_k$ of $k$-rough numbers, where $k$ is a prime number. Given $S_k$ for
any prime $k$, the set $S_{\kpr}$, where ${\kpr}$ is the smallest prime
greater than $k$, can be constructed easily. From this it follows that
PiFormulas occur in disjoint families. In any family, there is a 
first member, a series  of a
least prime number $k_{min}$ , which must be summed from first principles.
Then from the
series of some $k > k_{min}$ already summed, the series for ${\kpr}$ 
can be summed by a simple algebraic procedure. A good number of Pi Formulas,
belonging to a variety of families, and
giving results believed to be new, are presented here.

\end{abstract}

$\quad$

$\quad$

\nit {\bf Keywords}: rough number, infinite series, prime number, 
modulo multiplication group, totient function, contour integral

$\quad$
 
\nit {\bf Mathematical Subject Classification}: 11A67, 11B50, 40C10, 40C15
$\quad$

$\quad$


\section{Introduction}

This article presents a number of Pi Formulas believed to be new here, and
develops the view that Pi Formulas occur in families defined by 
association with the family of rough numbers.

Begin by quoting four Pi Formulas
\begin{align} 1-\frac{1}{3}+\frac{1}{5}-\frac{1}{7}+ \dots =&
\int^1_0 \frac{dx}{1+x^2} = \pi/4 \label{aa1} \\
1+\frac{1}{5}-\frac{1}{7}-\frac{1}{11}+ \frac{1}{13}+\frac{1}{17}-\frac{1}{19}
-\frac{1}{23}+\frac{1}{25}+\dots =&
\int^1_0 \frac{1+x^4}{1+x^6}dx = \pi/3\label{aa2} \end{align}
\begin{align} 
 & (1-\frac{1}{7}-\frac{1}{11}+\frac{1}{13}+\frac{1}{17}
-\frac{1}{19}-\frac{1}{23}+\frac{1}{29})- \nonumber \\
& (\frac{1}{31}-\frac{1}{37}-\frac{1}{41}+\frac{1}{43}+\frac{1}{47}
-\frac{1}{49}-\frac{1}{53}+\frac{1}{59}) + \dots   \nonumber \\
&= \int^1_0 \frac{(1-x^6)(1-x^{10})(1+x^{12})}{1+x^{30}} dx = \frac{4\pi}{15}
\label{aa3} \end{align}
\be\label{aa4}
1-\frac{1}{5}+\frac{1}{7}-\frac{1}{11}+\frac{1}{13}-\frac{1}{17}+\frac{1}{19}
-\frac{1}{23}+
\dots = \int^1_0
\frac{1-x^4}{1-x^6} dx= \frac{\pi\sqrt{3}}{6}
 \e\nit

Eq. (\ref{aa1}) is the famous result of Gregory (1671) and Leibniz (1674). The
second is found in \cite[eq. \ (36)]{eww3}. The third one is proved here in the
belief that it is a new result, as are various further results given below,
which differ from (\ref{aa3}) only in the distribution of signs of
the same set of fractions. Eq. (\ref{aa4}) is (85) in \cite{jol}.

Two good and extensive sources of Pi Formulas are the book of Jolley
\cite{jol} and \cite{eww3}. Eq. (\ref{aa2}) is not given in \cite{jol},
although the related result (\ref{aa4}), is, and 
\be\label{aa5}  1-\frac{1}{5}-\frac{1}{7}+\frac{1}{11}+
\frac{1}{13}-\frac{1}{17}-\frac{1}{19}+\frac{1}{23}+ \dots= \int^1_0
\frac{1-x^4}{1+x^6} dx=
\frac{1}{\sqrt{3}}
\log(2+\sqrt{3}), \e\nit
not of course a Pi Formula, appears as (84) there. Note also  exercise (528) in
\cite{hsc}.

The material presented in this article is organised as follows. Since rough
numbers and to a lesser extent modulo multiplication groups play an important
role in the formalism developd here, Section two provides background
information and references for these topics. Section three develops the view
that Pi Formulas occur in families in which passage from one member of the
family to its successor follows a uniquely defined path. Section four
discusses the evaluation of the integrals arising in Pi Formulas
such as those quoted above, and also below. Elementary methods,
contour integration, and an elegant algebraic procedure,  described in this
section, are all essential to the systematic accumulation of Pi Formula
results. The last of these methods is used to
extend the $\Sigma_3$ family of (\ref{aa1}) to an $S_{11}$ fourth member, and
potentially further. Section five, entitled More Pi Formulas, covers a range of
further examples which illustrate themes developed in previous sections,
but no systematic survey of the scene is attempted. Section six considers
briefly further developments.

\section{Rough numbers and Modulo Multiplication Groups (MMG)}

The $k$-rough numbers \cite{eww1}, defined for all prime numbers $k \geq 2$,
are the elements of the set $S_k$, which contains exactly once each positive
integer not divisible by any prime less than $k$.
Explicitly, showing also reference numbers in \cite{oeis},
\begin{align}S_3 = & \{ 1,3,5,7,9,11,13,15,17 \dots \},\quad
\underline{A005408}, \nonumber \\
S_5 = & \{ 1,5,7,11,13,17,19, \dots \} \quad \underline{A007310}, \nonumber \\
S_7 = & \{ 1,7,11,13,17,19,23, \dots \} \quad \underline{A007775} \nonumber \\
S_{11} = & \{ 1,11,13,17,19,23, \dots \} \quad \underline{A008364} \label{aa6}.
\end{align} \nit
The sequences $S_{13}$ and $S_{17}$ have references numbers
$\underline{A008365}$
and $\underline{A008366}$ in \cite{oeis}.
Note that $25,35,55,65,... $ belong to $S_5$, $49,77,119,... $ to $S_7$,
and $121,143, \dots $ to $S_{11}$, and so on. Also $S_2=\N$, so that all
natural numbers are $2$-rough, and  $S_1$ is undefined.

There is no entry rough number in the Mathematical Subject Classification,
although they feature as noted in \cite{oeis}.

It is indicated in \cite{eww1}, that the $k$-rough numbers were originally
introduced
in \cite{sf1}, \cite[Sec. 2.21]{sf2}. However it is clear that Euler was
well aware of rough number patterns, but, they are dealt briefly with
in passing, in the context of
studies of the harmonic series, a study in which divergent analogues of
(\ref{aa1} -- \ref{aa3}) appear. See \cite[p. 58]{wd}. It is furthermore
implicit there that the $S_k$ are easily constructed in order of
increasing $k$ by means
of the result
\be\label{aa65} S_{\kpr} =S_k-kS_k, \e\nit
where $\kpr$ is the smallest prime greater than $k$.

Turning to modulo multiplication groups (MMG), in relation to the $S_k$, define,
for each prime $k$, $P_k$ to be the product of all primes less than $k$;
in particular
\be\label{aa7} P_3=2, \quad P_5=6, \quad P_7=30, \quad P_{11}=210, \e\nit
and $P_{\kpr}=kP_k$. 
Now, if the elements of $S_k$ are reduced modulo $P_k$, the distinct integers
obtained provide the elements of a group $G_k=M(P_k)$ of order $g_k=\phi(P_k)$,
where $\phi(n)$ is the Euler totient function \cite{ds}, and the group
multiplication law is multiplication modulo $P_k$. The group $M(n)$ is referred
to \cite[p. 61]{ds},\cite{eww2} as a modulo multiplication group (MMG).
For example
\be\label{aa8}  G_3=M(2)=\{1 \},\; G_5=M(6)=\{1,5 \}, \; G_7=M(30)=
\{1,7,11,13,17,19,23,29 \}
\e\nit
of orders $g_3=1,\; g_5=2,\; g_7=8$. Also $G_{11}=M(210)$ is of order
$g_{11}=48$, and $G_{13}=M(2310)$ is of order $g_{13}= 480$. Note also
\be\label{aa9} g_{\kpr}=\phi(k)g_k=(k-1)g_k, \e\nit
as $k$ is prime.

The first $g_k$ element of each $S_k$ are the elements of
$G_k$. The remainder may be written down in successive sets of size $g_k$
obtained by adding $P_k,2P_k,\dots$ to elements of the first sets. For example,
\be\label{aa10} 1,5; 7,11; 13,17; 19,25; 29,31;35,37; \dots \e\nit

In all of (\ref{aa1} -- \ref{aa5}), in the denominators of the fractions summed
over, there is found the complete
set of elements of the relevant $S_k$. So (\ref{aa1}) is an $S_3$ result. The
four that follow in Section one are $S_5,S_7,S_5,S_5$ results.

\section{Families of $S_k$ results}

It is argued here that (\ref{aa1} -- \ref{aa3}) provide the $S_k$ results of
the $\Sigma_3$ family of such results, for $k=3,5,7$. Results for higher $k$
are obtained further on. Also (\ref{aa4}) is the first or $k=5$ result  of a
family $\Sigma_5$ of results, independent of the $\Sigma_3$ family.
The second member of this $\Sigma_5 $ family is found as (\ref{kk1})
in Section four.

To justify the view being advanced  and account for the statements just made,
begin by making some definitions:

$\quad$ i) the functions $f_k(x), \quad k=3,5, \dots$, which are the
integrands of the integrals which yield the series of the
family $\Sigma_3$, and which are given, for $k=3,5,7$ in 
(\ref{aa1} -- \ref{aa3})

$\quad$ ii) the functions $g_k(x), \quad k=3,5, \dots$, which play a similar
role for the series $\Sigma_5$, with $g_5(k)$ given by (\ref{aa4})

$\quad$ iii) the functions $h_k(x), \quad k=3,5, \dots$, associated with
(\ref{aa5}). See Sec. 5, last paragraph.

Expansions of the $f_k(x)$ give power series which,
upon integration, yield the series on the left side of
(\ref{aa1} -- \ref{aa3}). As noted above the fractions thereby found contain
all the elements of the corresponding $S_k$. When the terms of (\ref{aa1})
containing integers divisible by three are excised from it, only the
terms of (\ref{aa2}) remain, with the signs inherited from (\ref{aa1}).
Similarly, when the terms of (\ref{aa2})
containing integers divisible by five are excised from it, only the
terms of (\ref{aa3}) remain, with the signs inherited from (\ref{aa2}),
and hence from (\ref{aa1}). These statements indicate the sense in which
$S_3, S_5,S_7$  are the members of the $\Sigma_3$ family with ancestor $S_3$,
and in which $S_5$ is the child of its $S_3$ parent, and $S_7$ is
the grandchild of $S_3$. The continuation of the family to realise
further descendents follows a clear and unique path. As noted, (\ref{aa4})
is not a descendent of (\ref{aa3}), but is the parent of 
(\ref{kk1}).

Some further comments address tne complexity of the results of the
 $S_{11}$ or $S_{13}$ members of the family $\Sigma_3$. The first two terms
 of (\ref{aa2}), $\frac{1}{1}+\frac{1}{5}$ contain the integers of $G_5$.
 Then the denominator $(1-x^6)=(1-x^{P_5})$ of $f_5(x)$ and integration ,
 generates the rest of the (\ref{aa2})
\be\label{aa11} (1+\frac{1}{5})-(\frac{1}{7}+\frac{1}{11})+
(\frac{1}{11}+\frac{1}{13})-\dots \e\nit
Similarly, the first bracketed set of $g_8$ terms of (\ref{aa3}) contain the 
integers of $G_8$, and the denominator  $(1-x^{30})=(1-x^{P_7})$ of $f_7(x)$,
and integration, generates the rest of (\ref{aa3}). For the $k=11$ member
of the $\Sigma_3$ family, $g_{11}=48$ and the denominator function is 
$(1-x^{210})=(1-x^{P_{11}})$. This case is addressed in full below.

\section{Evaluation of integrals}

Turn next to the evaluation of integrals in  (\ref{aa1} -- \ref{aa4}), and
thereby to completion of the proofs of these four Pi Formulas. The integrals in
(\ref{aa1}) and (\ref{aa2}) can be evaluated by the standard procedures of
the integral calculus. But (\ref{aa3}) cannot be treated in this way. 
However conversion of
the integrals involved
in all four cases into contour integrals around a closed semicircular
contour $C$ in the upper half complex plane allows residue theory to be
employed. This offers a direct path to verifying the truth of (\ref{aa2}) and
(\ref{aa3}) and others displayed below. This is illustrated now for (\ref{aa2}).

To obtain a complex integral of the type required for (\ref{aa2}), note that
\be\label{bb1}
\int^1_0 \frac{1+x^4}{1+x^6} dx=\int^{\infty}_1 \frac{1+x^4}{1+x^6} dx
= \frac{1}{4} \int^{\infty}_{-\infty} \frac{1+x^4}{1+x^6} dx \e\nit
follows use of a change $y=1/x$ of variables for the first step. This leads to
the complex integral
\be\label{bb2} \oint_C \frac{1+z^4}{1+z^6} dz=
\oint_C \frac{g(z)}{h(z)} dz \e\nit
where $h(z)=1+z^6$ has simple poles within $C$ for $z=\exp(\fract{i\pi}{6}),
\exp(\fract{i\pi}{2}),
\exp(\fract{5i\pi}{6})$. If $z=e^{i\theta}$ is of any one of these, its residue
is given \cite{wik} by
\be\label{bb3} \frac{g(z)}{h^{\prime}(z)}= \frac{2\cos(2\theta)}{6z^3}, \e\nit
and the correctness of (\ref{aa2}) is confirmed, much more easily than by
more elementary methods.

Now turn to (\ref{aa3}), using
\be\label{cc1} J=\oint_C \frac{(1-z^6)(1-z^{10})(1+z^{12})}{(1+z^{30})}dz.
\e\nit
For $J$ there are $15$ poles inside the contour $C$ for
\be\label{cc2} z_r= \exp(\fract{(2r-1)i\pi}{30}), \quad r=1,2, \dots, 15.
\e\nit

If $z=e^{i\theta}$ is any one of these, its residue is
\be\label{cc3} \frac{-8\sin(3\theta)\sin(5\theta)\cos(6\theta)}{30z^{15}} \e\nit
Also $z^{15}= \pm i$ for all these poles. Proof that (\ref{aa3}) is correct
can now be completed.

The fourth Pi Formula in the $\Sigma_3$ family of (\ref{aa1} -- \ref{aa4})
is considered next. For the series on its left hand side  the $g_{11}=48$
terms of its first bracketed set of terms follows from the series of (\ref{aa3})
by excising all fractions with a denominator divisible by $7$, and extended  to
cover all of $S_{11}$ by use of the denominator $(1+x^{P_{11}})=(1+x^{210})$
and integration. The result can be expressed in the form
\be\label{hh0} K= \int^1_0 f_{11}(x)dx=\int^1_0 \frac{\psi_{11}(x)}
{1+x^{210}} dx \e\nit
where $\psi_{11}(x)$ is a polynomial degree $208$ with $48$ terms.
It seems to be futile to try to use contour integration to evaluate $K$,
since (a) no factorisation or other simplification of $\psi_{11}(x)$ has
been found, and (b) $C$ has $105$ poles within it. A new approach is called for.

A suitable approach emerges as follows. Begin by stating some results that
are easy to verify. First
\be\label{hh2} f_5(x)=f_3(x)+x^2 f_3(x^3), \e\nit
which gives,  upon using the change $y=x^3$ to treat the integral
of the second term, the result
\be\label{hh3} \int^1_0 f_5(x)dx=\int^1_0 f_3(x) dx (1+\frac{1}{3})=
\frac{\pi}{3}, \e\nit
in agreement with (\ref{aa2}). Second
\be\label{hh4} f_7(x)=f_5(x)-x^4 f_5(x^5), \e\nit
gives,  upon using the change $y=x^5$ to treat the integral
of the second term,
\be\label{hh5} \int^1_0 f_7(x)dx=\int^1_0 f_5(x) dx (1-\frac{1}{5})=
\frac{4\pi}{15}, \e\nit
in agreement with (\ref{aa3}). This prompts consideration of
\be\label{hh6}
f_{11}(x)=f_7(x)+x^6 f_7(x^7). \e\nit
For the left side of (\ref{hh6}), see (\ref{hh0}). Since $f_7(x)$ is known,
the right hand side of (\ref{hh6}) can be calculated using (\ref{hh4}).
Agreement with the
findings for the left hand side has been shown to hold, certainly as far as
the first $96$ terms, that is, the first two bracketed sets of $48$ terms,
are concerned.
Hence, for the fourth or $S_{11}$ member
of the $\Sigma_3$ family, the result
\be\label{hh7} \int^1_0 f_{11}(x) dx= \int^1_0 f_7(x) dx (1+\frac{1}{7})=
\frac{32\pi}{105} \e\nit
is advanced.

The approach just developed also gives rise to another result,
the $S_7$ descendent Pi Formula of (\ref{aa4}). It leads directly
to
\begin{align} 
 &(1+\frac{1}{7}-\frac{1}{11}+\frac{1}{13}-\frac{1}{17}
+\frac{1}{19}-\frac{1}{23}-\frac{1}{29})+ \nonumber \\
& (\frac{1}{31}+\frac{1}{37}-\frac{1}{41}+\frac{1}{43}-\frac{1}{47}
+\frac{1}{49}-\frac{1}{53}-\frac{1}{59}) + \dots   \nonumber \\
&= \int^1_0 g_7(x)dx.
\label{kk1} \end{align} \nit
where
\be\label{kk2} g_7(x)=\frac{(1+x^6)(1-x^{10})(1+x^{12})}{(1-x^{30})}. \e\nit
Hence, reading $g_5(x)=(1-x^4)/(1-x^6)$ off (\ref{aa4}),
\be\label{kk3} g_7(x)=g_g(x)+x^4g_5(x^5) \e\nit
follows, and the series of (\ref{kk1}) sums to
\be\label{kk4} \int^1_0 g_5(x)dx (1+\frac{1}{5})=
\frac{\sqrt{3}\pi}{5}. \e\nit
Contour integration confirms this.

There is no practical obstacle to treating the $S_{11}$ grandchild of
(\ref{aa4}).

A further point to note. Since (\ref{aa65}) deals with sets of
positive integers, the sign of its second term is always negative. Pi Formulas,
on the other hand, contain series with well-defined patterns of signs, which,
in formulas such as (\ref{hh2}), (\ref{hh4}) and (\ref{kk3}) require a
case by case check on the
correct sign to use for their second terms.

\section{More Pi Formulas}

Before turning attention to the discovery of Pi Formulas beyond the first
few or perhaps only the first member of any family of such formulas using the
fast algebraic approach, it should be noted that elementary methods or
contour integation provide the only way to get any series started.

For the discussion of $S_7$ results to follow, some suitable notation is
called for to avoid writing out results like those of (\ref{aa3}) and
(\ref{kk1}) in their displayed extended forms. All the results to be
stated involve the same sets of eight fractions in their successive
bracketed sets. But the distribution of the signs differ from case to case.
Clearly the series of (\ref{aa3}) and (\ref{kk1}) are fully specified by the
sign patterns of their fractions via
\begin{align} & (+--+ +--+)- \nonumber \\
&  (++-+ -+--)+ \label{kk5} \end{align} \nit
where the right most signs indicate results involving $(1 \pm x_{30})=
(1 \pm x_{P_7})$ respectively.

This enables a reasonably economical presentation of the following set of eight
$S_7$ results

\begin{align}
&  (+ + + +  + + + +)-  &  \frac{(1+x^6)(1+x^{10})(1+x^{12})}{1+x^{30})}
      &  \qquad \frac{4\sqrt{3}\pi}{15} \cos(\pi/10) \label{ss0}     \\
&  (+ - - +  + - - +)-  &  \frac{(1+x^6)(1-x^{10})(1-x^{12})}{1+x^{30})}
      &  \qquad \frac{4\pi}{15}            \label{ss1}     \\
&  (+ - + -  - + - +)-  &  \frac{(1-x^6)(1+x^{10})(1-x^{12})}{(1+x^{30})}
      &  \qquad \frac{4\sqrt{3}\pi}{15} \sin(\pi/5) \label{ss2}      \\
&  (+ + - -  - - + +)-  &  \frac{(1+x^6)(1-x^{10})(1-x^{12})}{(1+x^{30})}
      &  \qquad \frac{2\sqrt{5}\pi}{15}             \label{ss3}        \\
&  (+ - - -  + + + -)+  &  \frac{(1+x^6)(1-x^{10})(1+x^{12})}{(1-x^{30})}
      &  \qquad \frac{\pi}{\sqrt{15}}               \label{ss4}        \\
&  (+ + - +  - + - -)+  &  \frac{(1+x^6)(1-x^{10})(1+x^{12})}{(1-x^{30})}
      &  \qquad \frac{\sqrt{3}\pi}{5}               \label{ss5}        \\
&  (+ - + +  - + - -)+  &  \frac{(1-x^6)(1-x^{10})(1+x^{12})}{(1-x^{30})}
      &  \qquad \frac{\pi}{15} \sqrt{25-2\sqrt{5}}   \label{ss6}     \\
&  (+ + + -  + - - -)+  &  \frac{(1+x^6)(1+x^{10})(1-x^{12})}{(1-x^{30})}
      &  \qquad \frac{\pi}{15} \sqrt{25+2\sqrt{5}}   \label{ss7}
\end{align}
      \nit
The sums to the rightmost end of each line are those of the series specified
by the sign patterns at the leftmost end, using the integrands shown. All
of these results have been obtained by contour integration. Eq. (\ref{ss1})
is identical to (\ref{aa3}), and  
Eq. (\ref{ss5}) to the result (\ref{kk4}) obtained above for the $S_7$ child of
(\ref{aa4}).

Informallly it might be suggested that some of the results above merit the
label smooth. In order of smoothness therefore perhaps place
(\ref{aa3}) first, and follow with (\ref{ss4}), (\ref{ss5}) and (\ref{ss3}).

In all the examples displayed on the $S_7$ scene so far, all the
integrands used contain
an even number of minus signs, else the first step of the analogue of
(\ref{bb1}) fails. Eq. (\ref{aa5}), for example, can be evaluated by
elementary means but not by contour integration. 
Eq. (\ref{aa5}), however, is
the first member of a family of series, not of course providing Pi Formulas,
in which the method used
above to pass to further members of the family can nevertheless be employed.
It yields an $S_7$ series with sign pattern
\be\label{hh12} (+ - + +  - - + -)- \e\nit
and integral
\be\label{hh13} \int^1_0 h_7(x)dx, \;\; h_7(x)=\frac{(1-x^6)(1+x^{10})
(1+x^{12})}{(1+x^{30})}, \e\nit
which cannot be evaluated using contour integration.
Reading the integrand $h_5(x)=(1-x^4)/(1+x^6)$ off (\ref{aa5}), the result
\be\label{hh14} h_7(x)=h_5(x)+x^4 h_5(x^5), \e\nit
follows, so that the $S_7$ series under consideration sums to
\be\label{hh15} \int^1_0 h_7(x)dx= \int^1_0 h_5(x)dx (1+\frac{1}{5})=
\frac{2\sqrt{3}}{5} \log(2+\sqrt{3}).\e\nit
The last paragraph provides the last of the nine $S_7$ series of known sum to
hand.

\section{Further developments}

The PiFormulas given above display series defined by the simplest integrals
appropriate to their $S_k$ types. Ues of more complicated ones opens avenues
to a range of further PiFormulas. A few illustrations follow.

First the $S_3$ result
\be\label{jj1}   1+\frac{1}{3}-\frac{1}{5}-\frac{1}{7}+\frac{1}{9}+\frac{1}{11}
-\frac{1}{13}-\frac{1}{15}+\frac{1}{17}+\frac{1}{19}- \dots= \int^1_0
\frac{1+x^2}{1+x^4} dx=
\frac{\sqrt{2}\pi}{4}
 \e\nit
Can this smooth result possibly be new?
Next the $S_5$ result
\begin{align} 
 & (1+\frac{1}{5}+\frac{1}{7}+\frac{1}{11}+\frac{1}{13}+\frac{1}{17}
+\frac{1}{19}+\frac{1}{23})- \nonumber \\
& (\frac{1}{25}+\frac{1}{29}+\frac{1}{31}+\frac{1}{35}+\frac{1}{37}+
\frac{1}{41}+\frac{1}{43}+\frac{1}{47})
+ \dots   \nonumber \\
&= \int^1_0 \frac{(1+x^4)(1+x^6)(1+x^{10})}{(1+x^{24})} dx =
\frac{\pi}{3}(1+\sqrt{2})^{1/2}
\label{jj2} \end{align}
Finally, both by contour integration and an algebraic method, for  the $S_5$
child of (\ref{jj1}) find
\be\label{jj3}  (1-\frac{1}{5}-\frac{1}{7}+\frac{1}{11})
-(\frac{1}{13}-\frac{1}{17}-\frac{1}{19}+\frac{1}{23})+ \dots= \int^1_0
\frac{(1-x^4)(1-x^6)}{(1+x^{12})} dx=
\frac{\sqrt{2}\pi}{6}, \e\nit
resembling, but of course distinct from (\ref{aa5}).

\end{document}